\documentclass[12pt,english]{article}
\def\IC{\relax\,\hbox{$\inbar\kern-.3em{\rm C}$}}
\newcounter{lemma}[section]
\newcounter{proof}[section]
\newcounter{claim}[section]
\newcounter{corollary}[section]
\newcounter{theorem}[section]
\newcounter{proposition}[section]
\newcounter{definition}[section]
\newcounter{example}[section]
\newcounter{problem}[section]
\newcounter{remark}[section]
\newcounter{question}[section]
\def\vflx#1#2{\mathrel {\buildrel \hbox{$#1$} \over {\hbox to #2em
{\rightarrowfill}}}}
  \def\flx#1{\vflx{#1}{2}}
  \def\flxb{\Big\downarrow}

\date{}
\title{A noncommutative version of the Banach-Stone Theorem
}

\author{B. BOUALI \\University of Mohammed premier \\OUJDA, MOROCCO }

\begin{document}
\maketitle
\abstract{In this paper, we extend the Banach-Stone theorem to the non commutative case, i.e, we prove that the structure of  the liminal $C^{*}$-algebras $\cal A$ determines the topology of its primitive ideal space. 
}
\vspace{10mm}\\
A. M.S 2000 subject classification  : 46H05, 46H10, 46H15
\vspace{5mm}
\section{Introduction}

Let $X$ be a Banach space and let $C(S,X)$ (C(S)) denote the space of $X$-valued (Scalar-valued) continuous functions on a compact hausdorff space $S$ (endowed with the sup-norm). The classical Banach-Stone theorem states that the existance of an isometric isomorphism from $C(S)$ onto $C(S)$ implies that $S$ and $S'$ are homeomorphic.  There exists a variety of results in the literature linking the topological structure of a topological space $X$ with algebraic or topological-algebraic structures of $C(X)$, the set of all continuous real functions on $X$.  Further results along this line were obtained by Hewitt [1] and Shirota [2]. They proved respectively that, for a realcompact space $X$, the topology of $X$ is determined by the ring structure of $C(X)$ and by  the lattice structure of $C(X)$. Moreover, Shirota proved  in [2] that the lattices $UC(X)$ and $UB(X)$ determine the topology of a complete metric space $X$, where $UC(X)$ denotes the family of all uniformly continuous real functions on $X$, and $UB(X)$ denotes the subfamily of all bounded functions in $UC(X)$. Moreover Behrends [3] proved that if the centralizers (for the definition see also [3]) of $X$ and $Y$ are one-dimensional then the existance of an isometric isomorphism between $C(S,X)$ and  $C(S',X)$ implies that $S$ and $S'$ are homeomorphic. Cambern [4] proved that if $X$ is finite-dimensional Hilbert space and if $\Psi$ is an isomorphism of $C(S,X)$ onto $C(S',X)$ with $\mid\mid \Psi \Psi^{-1} \mid\mid < \sqrt{2}$ then  $S$ and $S'$ are homeomorphic. In [5] Jarosz proved that there is an isometric isomorphism between $C(S,X)$ and  $C(S',X)$ with a small bound iff $S$ and $S'$ are homeomorphic.

In the last few years there has been interest in the connection between the uniformity of a metric space $X$ and some further structures over $UC(X)$ and $UB(X)$. Thus, Araujo and Font in [6], using some results by Lacruz and Llavona
[7], proved that the metric linear structure of $UB(X)$ endowed with the $sup$-norm determines the uniformity of $X$, in
the case that $X$ is the unit ball of a Banach space. This result has been extended to any complete metric space $X$ by Hern\'{a}ndez [8]. Garrido and Jaramillo in [9] proved that the uniformity of a complete metric space $X$ is indeed characterized not only bu $UB(X)$ but also $UC(X)$. 

In this note, considering  $\cal A$ be $C^*$-algebra and Prim($\cal A$) be the space of primitive ideals, we prove that the structure of  $\cal A$ determines the topology of Prim($\cal A$).

\subsection{The hull kernel topology} The topology on Prim($\cal A$) is given by means of a closure operation. Given any subset $W$ of Prim($\cal A$), the closure $\overline{W}$ of $W$ is by definition the set of all elements in Prim($\cal A$) containing $\cap W=\lbrace \cap I  ~: ~ I \in W\rbrace $, namely 
$$
 \overline{W}=\lbrace I\in Prim({\cal A}): I \supseteq \cap W \rbrace$$

It follows that the closure operation defines a topology on Prim($\cal A$) which called Jacobson topology or hull kernel topology [10].

\vspace{3mm}\begin{proposition}[12, Proposition 2.7]The space Prim($\cal A$) is a $T_0$-space, i.e. for any two distinct points of the space there is an open neighborhood of one of the points which does not contain the other.\end{proposition}

\vspace{3mm}\begin{proposition}[12, Proposition 2.7] If $\cal A$ is $C^{*}$-algebra, then Prim($\cal A$) is locally compact.If $\cal A$ has a unit, then Prim($\cal A$) is  compact. \end{proposition}

\vspace{3mm}\begin{remark} The set of ${\cal K}({\cal H})$ of all compact operators on the Hilbert space $\cal H$ is the largest two sided ideal in the $C^*$-algebra  ${\cal B}({\cal H})$ of all bounded operators.\end{remark}

\vspace{3mm}\begin{definition}
A $C^*$-algebra $\cal A$ is said to be liminal if for every irreducible representation $(\pi, {\cal H})$ of $\cal A$ one has that $\pi(\cal A)= {\cal K}({\cal H})$
\end{definition}
\vspace{3mm}

So, the algebra $\cal A$ is liminal if it is mapped to the algebra of compact operators under any irreducible representation. Furthermore if $\cal A$ is a liminal algebra, then one can  prove that each primitive ideal of $\cal A$ is automatically a maximal closed two-sided ideal. As a consequence, all points of Prim($\cal A$) are closed and Prim($\cal A$) is a $T_1$-space . In particular, every commutative  $C^*$-algebra is liminal [10].
    
\section{The main result}
In this section, we extend the Banach-Stone Theorem to liminal $C^*$-algebras, before we give some lemma. 

\begin{lemma}
Let $\cal A$  and $\cal B$ be liminal $C^*$-algebras and let  $\alpha$ be an isomorphism of $\cal A$ onto $\cal B$. If $I$ is a primitive ideal of  $\cal B$ , then $\alpha^{-1} (I)$ is a primitive ideal of  $\cal A$. 
\end{lemma}

\vspace{3mm}
{\bf Proof} 
It is clear that the kernel of $\pi \times \alpha$( representation of $\cal A$ is $\alpha^{-1}(I_\pi)$, where $I_\pi$ is a primitive ideal of $\cal B$.

Now, we prove that $\pi\times \alpha$ is an irreducible representation of $\cal A$. If, contrary there exists a $\Pi(\cal A)$-invariant subspace $K$ of Hilbert space  $H$ ($K\not = 0 , K\not = H$)($\pi\times\alpha(\cal A) K=K$), a sample calcule show that $K$ is $\pi(\cal B$)-invariant and  $\pi$ is not a irreducible representation of $\cal B$. This is a contradiction, and we conclued that $\alpha^{-1}(I_\pi)$ is a primitive ideal of $\cal A$.       
\vspace{3mm}

\begin{theorem}
Let $\cal A$  and $\cal B$ be liminal $C^*$-algebras and let  $\alpha$ be an isomorphism of $\cal A$ onto $\cal B$. If $I$ is a primitive ideal of  $\cal B$ , then $\alpha^{-1} (I)$ is a primitive ideal of  $\cal A$. The map $I \to \alpha^{-1} (I)$ is a homeomorphism of Prim($\cal B$) onto Prim($\cal A$). 
\end{theorem}

\vspace{3mm}
{\bf Proof} The primitive ideal of $\cal B$ is a maximal ideal [10 ,corollary 4.1.11.(ii)], Let $I_\pi$ be a maximal ideal of $\cal B$ for some $\pi \in \hat {\cal B}$. From a lemma 2.1, $\alpha^{-1}(I_\pi)$ is a primitive ideal of $\cal A$ (then a maximal ideal of  $\cal A$, there is a function h
$$ h: Prim({\cal B}) \mapsto Prim({\cal A})$$\\ such that $\alpha^{-1}(I_\pi)=I_{h(\pi)}$. 

Since we can replace $\alpha^{-1}$ by $\alpha$, it follows that $h$ is a bijection. we have induced isomorphism $\chi_\pi : {\cal B}/I_\pi \mapsto \bf C$ given by $\chi_\pi (a)=\pi(a)$  and  $\beta : {\cal B}/I_\pi \mapsto {\cal A}/I_{h^{-1}(\pi)}$ given by $\beta(a+I_\pi ) =\alpha (a) +I_{h^{-1}(\pi)}$. Therefore we get a commutative diagram:

$$
  \matrix{
  {\cal B}/I_{\pi} & \flx\beta & {\cal A}/I_{h^{-1}(\pi)}\cr\cr
  \chi_\pi \flxb && \flxb \chi_{h^{-1}(\pi)} \cr\cr
  {\bf C} & \flx{\gamma} & {\bf C}
  }
  $$
\vspace{5mm}\\
 and an induced  isomorphism $\gamma : {\bf C} \mapsto {\bf C}$ defined by $$\gamma (\pi(a))=h^{-1}(\pi)(\alpha(a))$$ Since $\gamma$ must be the identity maps , then  $\pi(a)=h^{-1}(\pi)(\alpha(a))$.

From [11], all opens set of Prim($\cal A$) are of the forms:
 $$U_I =\lbrace P\in Prim({\cal A}) : P   \not \supseteq  I\rbrace$$
Calculate $h^{-1}(U_I)$.\\
$$\begin{array}{lll}
h^{-1}(U_I)&=& \lbrace  ~~~~\pi~~~  : ker(\pi) \in Prim({\cal A}) \mbox{~~and~~} ker(\pi)   \not \supseteq  I\rbrace\\
&=& \lbrace  h^{-1}(\pi) : ker(\pi)\in Prim({\cal A}) \mbox{~~and~~} ker(\pi)   \not \supseteq  I\rbrace\\
&=& \lbrace  ~~~~\pi '~~~ : ker(\pi ')\in Prim({\cal A}) \mbox{~~and~~} ker( h(\pi'))   \not \supseteq  I\rbrace\\
&=& \lbrace  ~~~~\pi '~~~ : ker(\pi ')\in Prim({\cal A}) \mbox{~~and~~} \alpha^{-1}(I_\pi)   \not \supseteq  I\rbrace\\
&=& \lbrace  ~~~~\pi '~~~ : ker(\pi ')\in Prim({\cal A}) \mbox{~~and~~} ker( \pi')   \not \supseteq  \alpha(I)\rbrace\\
&=& U_{\alpha(I)}.
\end{array}$$

Then $h^{-1}(U_I)$ is an open set of Prim($\cal B$) and h is continuous. Replace $\alpha$ by $\alpha^{-1}$, it follows that $h^{-1}$ is continuous, then $h$ is a homeomorphism.
\vspace{3mm}

\begin{remark} We can observe that when $\cal A$ and $\cal B$ are commutative, we can take ${\cal A}=C_0(X)$, ${\cal B}=C_0(Y)$, Prim(${\cal A})=X$, Prim(${\cal B})=Y$ and the above reduces to the famous Banach-Stone Theorem.
\end{remark} 

\vspace{3mm}\begin{remark} The set of ${\cal K}({\cal H})$ of all compact operators on the Hilbert space $\cal H$ is a liminal $C^{*}$-algebra and from theorem 2.1 has the Banach-Stone property.\end{remark}
\vspace{3mm}

\begin{question} Has the structure of  all $C^{*}$-algebras $\cal A$ determines the topology of its primitive ideal space ? 
\end{question} 
\vspace{3mm}

{\bf Acknowledgements.} It is a pleasure to thank
Professors M. El Hodaibi  for several
interesting conversations concerning this note. I am also
grateful to E. Benrends for his lecture Notes in Mathematics [3], to M. Isabel Garrido for his paper (with Jesus A. Jaramillo)[6], and to K. Jarosz for his paper [5].

\hspace{.2cm}

{\small \noindent Bouchta BOUALI

\noindent Facult\'e des Sciences,Departement de Mat\'ematiques

\noindent Universit\'e Mohammed Premier

\noindent 60000 Oujda, Maroc

\noindent bbouali@sciences.univ-oujda.ac.ma

\end{document}